\documentclass[11pt]{article}
\usepackage{geometry}  
\usepackage[round]{natbib}
\usepackage{url}

\newtheorem{theorem}{Theorem}
\newtheorem{lemma}[theorem]{Lemma}

\newtheorem{corollary}[theorem]{Corollary}
\newtheorem{proposition}[theorem]{Proposition}
\newenvironment{proof}{\emph{Proof}}{\strut} 
 
\usepackage[small]{caption}
\usepackage[doublespacing]{setspace}
\usepackage[english]{babel}
\usepackage[utf8x]{inputenc}
\usepackage[T1]{fontenc}

\usepackage{amssymb, amsmath} 

\newcommand{\0}{\mathbf{0}}
\newcommand{\balpha}{\boldsymbol{\alpha}}
\newcommand{\bOmega}{\boldsymbol{\Omega}}
\newcommand{\bGamma}{\boldsymbol{\Gamma}}
\newcommand{\bA}{\mathbf{A}}
\newcommand{\ba}{\mathbf{a}}
\newcommand{\bD}{\mathbf{D}}

\newcommand{\bdelta}{\boldsymbol{\delta}}

\newcommand{\bI}{\mathbf{I}}

\newcommand{\bj}{\mathbf{j}}

\newcommand{\bLambda}{{\boldsymbol{\Lambda}}}
\newcommand{\bmu}{\boldsymbol{\mu}}
\newcommand{\bP}{\mathbf{P}}
\newcommand{\bPsi}{{\boldsymbol{\Psi}}}

\newcommand{\bSigma}{\boldsymbol{\Sigma}}

\newcommand{\btau}{{\boldsymbol{\tau}}}
\newcommand{\btheta}{\boldsymbol{\theta}}

\newcommand{\bp}{\mathbf{p}}
\newcommand{\bU}{\mathbf{U}}
\newcommand{\bV}{\mathbf{V}}

\newcommand{\bx}{\mathbf{x}}
\newcommand{\bxi}{{\boldsymbol{\xi}}}

\newcommand{\by}{\mathbf{y}}
\newcommand{\bY}{\mathbf{Y}}
\newcommand{\bz}{\mathbf{z}}
\newcommand{\bZ}{\mathbf{Z}}

\newcommand{\Prob}{\text{Pr}}
\newcommand{\R}{\mathbb{R}}
\newcommand{\bDelta}{\boldsymbol{\Delta}}
\newcommand{\bomega}{\boldsymbol{\omega}}
\newcommand{\bnu}{\boldsymbol{\nu}}
\newcommand{\bG}{\mathbf{G}}

\newcommand{\Real}{\mathbb{R}}
\newcommand{\T}{^{\top}}

\begin{document}

\thispagestyle{empty} \baselineskip=28pt \vskip 5mm
\begin{center} {\Huge{\bf On the non-identifiability of unified skew-normal distributions}
}       
\end{center}

\baselineskip=12pt \vskip 10mm

\begin{center}\large
Kesen Wang,\footnote[1]{\baselineskip=10pt Statistics Program,
King Abdullah University of Science and Technology,
Thuwal 23955-6900, Saudi Arabia\\
E-mail: kesen.wang@kaust.edu.sa, marc.genton@kaust.edu.sa\\
This research was supported by the
King Abdullah University of Science and Technology (KAUST)
} 
Reinaldo B. Arellano-Valle,\footnote[2]{\baselineskip=10pt 
Departamento de Estadística, Pontificia Universidad Católica de Chile, Santiago, Chile\\
E-mail: reivalle@mat.uc.cl
}
Adelchi Azzalini,\footnote[3]{\baselineskip=10pt 
Dipartimento di Scienze Statistiche, Università degli Studi di Padova, Padua, Italy\\
E-mail: azzalini@stat.unipd.it
}
\\ and Marc G.~Genton{$^1$}
\end{center}
\baselineskip=17pt \vskip 10mm \centerline{\today} \vskip 15mm

\abstract{\noindent We investigate the non-identifiability of the multivariate unified skew-normal distribution under permutation of its latent variables. We show that the non-identifiability issue also holds with other parameterizations and extends to the family of unified skew-elliptical distributions and more generally to selection distibutions. We provide several suggestions to make the unified skew-normal model identifiable and describe various sub-models that are identifiable.}

\vspace{4ex}

\noindent 
\emph{Keywords:~}Latent variable, Non-identifiability, Permutation, Selection distribution, Unified skew-elliptical distribution, Unified skew-normal distribution.
 
\clearpage

\section{The Skew-Normal Distribution and Some Extensions}
\subsection{The broad background} \label{intro-basics}

In the last 20--25 years, there has been a growing interest in a formulation for the construction of probability distributions of which the archetypal instance is the skew-normal  (SN) distribution, among many other constructions. A recent, relatively concise  account of the already vast literature on this theme is provided by \cite{azzalini:2022}. 

Since the present contribution deals with an extension of the multivariate skew-normal distribution and an associated issue which will be presented later on, we start by recalling the SN basic constituents as a preliminary step. For expository simplicity, we  consider initially its `normalized' case  where the location parameters are all $0$ and the scale parameters are all set to $1$.   
If $\bar\bOmega$ is a symmetric positive definite correlation matrix and $\phi_d(\bz; \bar\bOmega)$ denotes the ${\cal N}_{d}(\0,\bar\bOmega)$ density function evaluated at $\bz\in\Real^d$, then the SN density with $\0$ location is obtained by multiplication of $\phi_d(\bz;\bar\bOmega)$ and a modulation factor $2\,\Phi(\balpha\T\bz)$, involving the ${\cal N}(0,1)$ distribution function, $\Phi$, evaluated at a suitable point $\balpha\T\bz$. The coefficients $\balpha$ regulate the non-normality of the resulting density, and a null vector $\balpha$ produces the classical multivariate normal density.

There exist various stochastic constructions leading to a SN random variable, of which we sketch the one directly relevant for our subsequent development. Start from a ${\cal N}_{1+d}(\0,\bar\bOmega^*)$ variable with standardized components, of which the first one is denoted $U_0$ and the others are arranged in a vector $\bU_1$ 
such that $\mathrm{Var}(\bU_1)=\bar\bOmega$. Next, introduce a selection mechanism  $U_0>0$ and consider the conditional distribution of $\bU_1$ given that the latent variable $U_0$ is positive; this conditional distribution is of type SN whose parameter $\balpha$ is an appropriate function of $\mathrm{Cov}(U_0, \bU_1)$ and $\bar\bOmega$.

Some years after the multivariate SN distribution had appeared, various authors have proposed extensions where the modulation factor $2\,\Phi(\balpha\T\bz)$ is replaced by one involving the $m$-dimension normal distribution function, $\Phi_m$, for $m\ge1$. Specifically, an expression of type $\Phi_m(\bx;\bSigma)$ will refer to the distribution function of $\phi_m(\bx; \bSigma)$ at $\bx\in\Real^m$.
A first instance of this type is the closed skew-normal (CSN) family proposed by \citet{gonz:etal:2004-inSE}. The CSN construction starts from a $(m+d)$-dimensional normal variable, of which $m$ components are imposed to satisfy one-sided inequality constraints; these constraints generalize the condition $U_0>0$ appearing above. The resulting density function, now with arbitrary location and scale parameters, is 
\begin{equation}
    \phi_d(\by-\bmu;\bSigma)\: \frac{\Phi_m(\bD(\by-\bmu)-\bnu ;\bDelta)}{\Phi_m(-\bnu; \bDelta + \bD\bSigma\bD\T)},
    \qquad \by\in\Real^d,
    \label{csn-pdf}
\end{equation}
 where $\bmu \in \R^{d}$, $\boldsymbol{\nu} \in \R^m$,  $\bD \in \R^{m \times d}$, $\bSigma \in \R^{d \times d}$ and $\bDelta \in \R^{m \times m}$, with the condition that  $\bSigma$ and $\bDelta$ are symmetric positive definite matrices, written $\bSigma>0$ and $\bDelta>0$, while $\bD$ is arbitrary.
 Density (\ref{csn-pdf}) is parameterized as ${\cal CSN}_{d,m}(\bmu,\bSigma,\bD,\boldsymbol{\nu},\bDelta)$, provided $\bSigma>0$ and  $\bDelta>0$.

\cite{arellano2005fundamental} have examined two similar constructions,   denoted `Fundamental SN' and `CSN-2' which share with  (\ref{csn-pdf}) the presence of a $\Phi_m$ term depending on $\by$ at the numerator and a matching normalizing constant which also involves the function $\Phi_m$.
Yet another construction of similar type was generated by \cite{lise:lope:2003}, under the name `Hierarchical SN'.

\subsection{Summary facts of the SUN distribution and the present aim}
\label{intro-sun}

The formulations recalled in the second half of Section~\ref{intro-basics} are visibly very similar but their specific connections are not obvious. This exploration was the initial motivation of \cite{AA06}, leading to two conclusions: (i)~the various distributions are essentially coincident under an appropriate matching of their parameterizations; (ii)~in some form or another, they all are affected by over-parameterization, leading to non-identifiability of the model when one comes to their use for inferential purposes.
For instance, statement (ii) in the case of the CSN class means that, for any diagonal matrix $\bG$ having all diagonal elements positive, ${\cal CSN}_{d,m}(\bmu,\bSigma,\bD,\boldsymbol{\nu},\bDelta)$  coincides with ${\cal CSN}_{d,m}(\bmu,\bSigma, \bG\bD, \bG\boldsymbol{\nu},\bG\bDelta\bG)$; for the proof of this fact and analogous ones for other formulations; see \citet[Section 2.2]{AA06}. The source of the problem is, in all these constructions, that one can arbitrarily scale the latent normal component on which the underlying selection mechanism is applied.

To overcome the above-mentioned over-parameterization problem, \cite{AA06} put forward a formulation, denoted  unified skew-normal (SUN), which embraces the previous ones and prevents the arbitrary scaling of the variables involved in the selection mechanism. Specifically, start from 
\begin{align}
    \bU = \begin{pmatrix}\bU_0\\
              \bU_1\end{pmatrix} 
    \sim {\cal N}_{m+d}(\0,\bar\bOmega^*), \qquad 
    \bar\bOmega^* = \begin{pmatrix}
                \Bar{\bGamma} & \bDelta^\top \\
                \bDelta & \Bar{\bOmega}
    \end{pmatrix} \label{def}
\end{align}
where $\bar\bOmega^*$ is a correlation matrix. The selection mechanism leading to skewness of the distribution is now $\bZ = (\bU_1|\bU_0 + \btau > \0)$ for some vector of truncation parameters $\btau\in\Real^m$; the inequality sign in $\bU_0 + \btau > \0$ and similar expressions must be intended to hold component-wise.
Note that $\bDelta$ has different meanings in (\ref{def}) and in (\ref{csn-pdf}), having retained the symbols used in the original publications. 

Finally, consider the affine transformation $\bY=\bxi+\bomega\,\bZ$, where $\bomega$ is a diagonal matrix with $\bomega>0$, and define  $\bOmega=\bomega\Bar{\bOmega}\bomega$. Then, the density function of $\bY$ is
\begin{align}
    f(\by)=\phi_d(\by-\bxi;\bOmega)\:\frac{\Phi_m(\btau + \bDelta^\top\Bar{\bOmega}^{-1}\bomega^{-1}(\by-\bxi);\Bar{\bGamma}-\bDelta^\top\Bar{\bOmega}^{-1}\bDelta)}{\Phi_m(\btau;\Bar{\bGamma})},
    \qquad \by\in\Real^d,
    \label{sun-pdf}
\end{align}
which is said to be the ${\cal SUN}_{d,m}(\bxi,\bOmega,\bDelta,\btau,\Bar{\bGamma})$ density.
A number of formal properties of the SUN are provided by \cite{AA06}; numerous other properties have been developed by \cite{gupta2013some} and \cite{AA22}.

 In the SUN construction, it is the condition that $\bar \bOmega^*$ is a correlation matrix which prevents arbitrary scaling of the latent variable $\bU_0$.
 However, in the next section we show that this condition does not cope with all possible sources of  over-parameterization,  and hence with the non-identifiability connected to the selection mechanism. {This non-identifiability of the SUN is a fundamental problem that has been overlooked in the literature and needs to be resolved before any attempt of developing inferential procedures for the family of SUN distributions. }
 The final part of the paper deals with a range of ensuing aspects and developments.


\section{Non-Identifiability of the SUN under Permutation of Latent Variables}
\label{SUN_nonidentifiability}

Given their role in the following discussion, recall first some key facts about permutation matrices, as given for instance by \citet[pp.~25-26]{HJ87}.
A permutation matrix $\bP$ of order $m\geq 2$ is an  $m\times m$  matrix obtained by permuting the rows of the identity matrix of order $m$, denoted $\bI_m$, according to some permutation of the numbers $1$ to $m$; a permutation of the columns of $\bI_m$ is an equivalent building mechanism. Hence, $\bP$ has one entry equal to $1$ in each row and each column, and $0$s elsewhere.  Clearly, there are $m!$ such permutation matrices, for a given order $m$; denote this set by  ${\cal P}(m)$. Note that a permutation matrix is not necessarily symmetric. It can easily be verified that the product of a permutation matrix $\bP$ and its transpose $\bP\T$ equals to the identity matrix; hence permutation matrices are orthogonal and their determinant equals~$\pm 1$. 

These basic properties allow us to  equivalently define  ${\cal P}(m)$ in the alternative form ${\cal P}(m)=\{\bP \in \R^{m \times m}| \bP\bP^\top = \bP^\top\bP = \bI_m \text{ and } \bP\boldsymbol{1}_m = \boldsymbol{1}_m\}$, where $\boldsymbol{1}_m$ denotes the $m$-vector of all $1$'s. To see this fact, consider that each row $\bp_i^\top$ of $\bP$ must satisfy $\bp_i^\top\bp_i=1$, $\bp_i^\top\bp_j=0$ for $j\neq i$, and $\bp_i^\top \boldsymbol{1}_m=1$. Geometrically, these conditions describe the intersection of a unit hypersphere centered at the origin and a hyperplane orthogonal to $\boldsymbol{1}_m$, the solutions of which are the vectors $\bp_i$ with all entries $0$s except one~$1$, and are mutually orthogonal.

Let $\bP \in {\cal P}(m)$ and consider the selection representation $\bZ_\bP = (\bU_1|\bP\bU_0 + \bP\btau > \0)$ where the ingredients
\begin{align*}
    \bU_\bP = \begin{pmatrix}\bP\bU_0\\
    \bU_1\end{pmatrix} \sim {\cal N}_{m+d}(\0,\bar\bOmega^*_\bP), 
    \qquad \bar\bOmega^*_\bP = \begin{pmatrix}
        \bP\Bar{\bGamma}\bP^\top & \bP\bDelta^\top \\
        \bDelta\bP^\top& \Bar{\bOmega}
    \end{pmatrix} 
\end{align*}
have the same structure as in (\ref{def}).
Correspondingly, define  $\bDelta_\bP = \bDelta\bP^\top$, $\btau_\bP = \bP\btau$, and $\Bar{\bGamma}_\bP = \bP\Bar{\bGamma}\bP^\top$, which is still a correlation matrix.
Then, with the affine transformation $\bY_\bP=\bxi+\bomega\,\bZ_\bP$, it is clear that $\bY_\bP \sim {\cal SUN}_{d,m}(\bxi,\bOmega,\bDelta_\bP,\btau_\bP,\Bar{\bGamma}_\bP)$ with density function
\begin{align}
        f_\bP(\by)=\phi_d(\by-\bxi;\bOmega)\:\frac{\Phi_m(\btau_\bP + \bDelta_\bP^\top\Bar{\bOmega}^{-1}\bomega^{-1}(\by-\bxi);\Bar{\bGamma}_\bP-\bDelta_\bP^\top\Bar{\bOmega}^{-1}\bDelta_\bP)}{\Phi_m(\btau_\bP;\Bar{\bGamma}_\bP)},
        \qquad \by\in\Real^d.
        \label{sun_p-pdf}
\end{align}
In the denominator, it is immediate that $\Phi_m(\btau;\Bar{\bGamma}) = \Phi_m(\btau_\bP;\Bar{\bGamma}_\bP)$. 
Next, we consider the arguments of the $\Phi_m$ term in the numerator, namely
\begin{align*}
    \btau_\bP + \bDelta_\bP^\top\Bar{\bOmega}^{-1}\bomega^{-1}(\by-\bxi) & = \bP\btau +\bP\bDelta^\top\Bar{\bOmega}^{-1}\bomega^{-1}(\by-\bxi)=\bP\{\btau + \bDelta^\top\Bar{\bOmega}^{-1}\bomega^{-1}(\by-\bxi)\},\\
    \Bar{\bGamma}_\bP-\bDelta_\bP^\top\Bar{\bOmega}^{-1}\bDelta_\bP & = \bP\Bar{\bGamma}\bP^\top - \bP\bDelta^\top\Bar{\bOmega}^{-1}\bDelta\bP^\top=\bP(\Bar{\bGamma}-\bDelta^\top\Bar{\bOmega}^{-1}\bDelta)\bP^\top.
\end{align*}
Since we only apply permutations, we have that $\Phi_m(\btau + \bDelta^\top\Bar{\bOmega}^{-1}\bomega^{-1}(\by-\bxi);\Bar{\bGamma}-\bDelta^\top\Bar{\bOmega}^{-1}\bDelta)=\Phi_m(\btau_\bP + \bDelta_\bP^\top\Bar{\bOmega}^{-1}\bomega^{-1}(\by-\bxi);\Bar{\bGamma}_\bP-\bDelta_\bP^\top\Bar{\bOmega}^{-1}\bDelta_\bP)$. Therefore, $f(\by)=f_\bP(\by)$ for all $\by\in\Real^d$, which implies that ${\cal SUN}_{d,m}(\bxi,\bOmega,\bDelta,\btau,\Bar{\bGamma})\equiv{\cal SUN}_{d,m}(\bxi,\bOmega,\bDelta_\bP,\btau_\bP,\Bar{\bGamma}_\bP)$ for any permutation matrix $\bP$, confirming the non-identifiability claim. 

The key point of the above argument is the equality of the probabilities at the numerators (respectively, denominators) of (\ref{sun-pdf}) and (\ref{sun_p-pdf}).
The equivalence of these distribution functions, before and after application of a transformation associated to a permutation matrix, holds in general for any multivariate distribution, not only for probabilities associated to a Gaussian distribution. In essence, the reason is that applying a permutation matrix simply exchanges the labels of the components, but leaves otherwise unchanged the components. A more detailed argument in support of this statement is given next.

Our aim is to characterize the set of linear transformations of the latent variable $\bU_0$ which lead to the non-identifiability issue. To this end, we shall need to use Lemmas~\ref{lem1} and \ref{lem2} stated below. Given the general nature of Lemma~\ref{lem1}, its result may well exist in the literature, but we could not actually locate an instance; hence, a proof is provided here.
In this lemma and for the rest of the paper, the term `positive diagonal matrix' is used as a shorthand for `diagonal matrix with positive diagonal elements'.
Note incidentally that a positive diagonal matrix is a positive definite matrix, while in general the same does not hold true for a positive matrix.

\begin{lemma}
Given a non-singular matrix $\bA=(a_{ij})\in\R^{m\times m}$,
the equivalence
\begin{equation}
\bx \leq \by \quad\Longleftrightarrow\quad \bA\bx \leq \bA\by \label{inequality}
\end{equation}
holds for all $\bx=(x_1,\ldots,x_m)\T$ and $\by=(y_1,\ldots,y_m)\T$ if and only if $\bA=\bD\bP$ where $\bD$ is a positive diagonal matrix and $\bP$ is a permutation matrix, $\bP\in {\cal P}(m)$.
\label{lem1}
\end{lemma}

\begin{proof} First, assume that $\bA=\mathrm{diag}(a_{11},\ldots,a_{mm})$ is a diagonal matrix with $a_{ii}>0$, $i=1,\ldots,m$, and prove that (\ref{inequality}) holds. Start by noticing that, if we apply the transformation $\bA$ on both sides of the inequality $\bx \leq \by$, this operation is only scaling the individual components $x_i$ and $y_i$ by $a_{ii}>0$, for $i=1,\dots,m$. Hence, $\bx \leq \by \implies \bA\bx \leq \bA\by$. Next, suppose that $\bx_\bA \leq \by_\bA$, where $\bx_\bA = \bA\bx$ and $\by_\bA = \bA\by$, and notice that $\bA^{-1}=\mathrm{diag}(a_{11}^{-1},\ldots,a_{mm}^{-1})$. Hence, if we apply $\bA^{-1}$ on both sides of the inequality $\bx_\bA \leq \by_\bA$, this operation is still just scaling the individual components $(\bx_{\bA}) _i$ and $(\by_{\bA})_i$ by $a_{ii}^{-1}>0$. Consequently, we have that $\bx_\bA \leq \by_\bA \implies \bx \leq \by$.

Second, assume that $\bA\in {\cal P}(m)$ is a permutation matrix and prove that (\ref{inequality}) holds. The statement follows immediately by considering that $\bA$ is only changing the order of the individual components, simultaneously for $\bx$ and $\by$, and it can easily be reversed by applying the inverse transformation, associated to $\bA^\top$, which is another permutation matrix.

Third, assume that $\bA=\bD\bP$ where $\bD$ is a positive diagonal matrix and $\bP$ is a permutation matrix, $\bP\in {\cal P}(m)$ and prove that (\ref{inequality}) holds. If $\bP\in {\cal P}(m)$ then from the second paragraph above we have $\bx \leq \by \Longleftrightarrow \bP\bx \leq \bP\by$. Denote $\bx_\bP = \bP\bx$ and $\by_\bP = \bP\by$. If $\bD$ is a positive diagonal matrix then from the first paragraph above we have $\bx_\bP \leq \by_\bP \Longleftrightarrow \bD\bx_\bP \leq \bD\by_\bP$. Therefore, $\bx \leq \by \Longleftrightarrow \bx_\bP \leq \by_\bP \Longleftrightarrow \bD\bx_\bP \leq \bD\by_\bP \Longleftrightarrow \bD\bP\bx \leq \bD\bP\by$.

Reciprocally, assume now that (\ref{inequality}) holds, and prove $\bA=(a_{ij})$ must be of the form $\bA=\bD\bP$ where $\bD$ is a positive diagonal matrix and $\bP$ is a permutation matrix, $\bP\in {\cal P}(m)$. By the first implication of (\ref{inequality}), we have that $\bx \leq \by \implies \bA\bx \leq \bA\by$.  We must show that the latter inequality  holds for any $\bx$ and $\by$ in $\R^m$ if and only if $a_{ij} \geq 0$ for all $i,j=1,\dots,m$. The `if' part is obvious; the `only if' part follows from the following argument \emph{ad absurdum}. Assume that one or more  elements of $\bA$ are negative, $a_{ij}<0$ say, for some fixed row index $i$ and all $j$ in a certain non-empty set $J$ of indices, $J\subseteq \{1,\dots,m\}$. Now choose $\bx=\0$, the $m$-dimensional null vector, and $\by$ with elements $y_k=1$ if $k\in J$, and $y_k=0$ otherwise. With these choices of $\bx$ and $\by$, $(\bA\bx)_i=0$ while $(\bA\by)_i<0$, which violates the condition $\bA\bx \le \bA\by$. Since we have reached a contradiction,  we conclude that all the $a_{ij}$ elements have to be non-negative. 
Similarly, the condition that $\bx_\bA \leq \by_\bA \implies \bx \leq \by$ for any $\bx$ and $\by$ in $\R^m$ implies that $a^{ij} \geq 0$, for all $i,j=1,\dots,m$, where $\bA^{-1}=(a^{ij})$. According to \cite{matrixinverse}, a non-singular non-negative matrix $\bA$ whose inverse $\bA^{-1}$ is also non-negative can only be the product of a positive diagonal matrix and a permutation matrix.
\end{proof}

\begin{lemma}
Let $\bY\in \R^m$ be an absolutely continuous random vector with distribution function $F_\bY(\cdot)$,   and let $\bA=(a_{ij})\in\R^{m\times m}$ be a non-singular matrix. Define $\bY_\bA=\bA\bY$ and denote its distribution function by $F_{\bY_\bA}(\cdot)$. Then, for all $\by\in\R^m$, $F_\bY(\by)=F_{\bY_\bA}(\by_\bA)$  where $\by_\bA=\bA\by$ if and only if $\bA=\bD\bP$ where $\bD$ is a positive diagonal matrix and $\bP$ is a permutation matrix, $\bP\in {\cal P}(m)$.
\label{lem2}
\end{lemma}

\begin{proof} From Lemma~\ref{lem1}, it follows that:
$$F_\bY(\by)=\Prob(\bY \leq \by)=\Prob(\bA\bY \leq \bA\by)=\Prob(\bY_{\bA} \leq \by_\bA)=F_{\bY_\bA}(\by_\bA)$$
if and only if $\bA=\bD\bP$ where $\bD$ is a positive diagonal matrix and $\bP$ is a permutation matrix, $\bP\in {\cal P}(m)$.
\end{proof}

We are now ready to state our main result.

\begin{proposition}
Let $\bY\sim {\cal SUN}_{d,m}(\bxi,\bOmega,\bDelta,\btau,\Bar{\bGamma})$ and $\bY_\bP \sim {\cal SUN}_{d,m}(\bxi,\bOmega,\bDelta_\bP,\btau_\bP,\Bar{\bGamma}_\bP)$ with $\bDelta_\bP = \bDelta\bP^\top$, $\btau_\bP = \bP\btau$, and $\Bar{\bGamma}_\bP = \bP\Bar{\bGamma}\bP^\top$. Then the two distributions coincide, i.e. ${\cal SUN}_{d,m}(\bxi,\bOmega,\bDelta,\btau,\Bar{\bGamma})\equiv{\cal SUN}_{d,m}(\bxi,\bOmega,\bDelta_\bP,\btau_\bP,\Bar{\bGamma}_\bP)$, if and only if $\bP\in {\cal P}(m)$.
\label{prop1}
\end{proposition}

\begin{proof} The `if' part was already proved at the beginning of this section. The `only if' part follows from Lemma~\ref{lem2} and the restriction that $\Bar{\bGamma}_\bP$ must be a correlation matrix, hence forcing $\bD=\bI_m$ in Lemma~\ref{lem2}.
\end{proof}

\begin{corollary}
The SUN family of distributions (\ref{sun-pdf}) is non-identifiable if $m>1$.
\end{corollary}
\begin{proof}
This follows immediately from Proposition~\ref{prop1}  because, for any member of the SUN class, the members associated to all possible choices of $\bP\in {\cal P}(m)$ have the same distribution.
\end{proof}

We draw attention to the fact that the problem exists only for $m>1$, because no permutation is possible for $m=1$. Since the original SN distribution corresponds to the SUN with $m=1$ and $\tau=0$,  it is clear from this problem too.

\begin{corollary}
Let $\bY\sim {\cal SUN}_{d,m}(\bxi,\bOmega,\bDelta,\btau,\Bar{\bGamma})$ and $\bY_\bP \sim {\cal SUN}_{d,m}(\bxi,\bOmega,\bDelta_\bP,\btau_\bP,\Bar{\bGamma}_\bP)$ with $\bDelta_\bP = \bDelta\bP^\top$, $\btau_\bP = \bP\btau$, and $\Bar{\bGamma}_\bP = \bP\Bar{\bGamma}\bP^\top$. Let $\bY=(\bY_1^\top,\bY_2^\top)^\top$, $\bY_\bP=(\bY_{\bP1}^\top,\bY_{\bP2}^\top)^\top$, with $\bY_1,\bY_{\bP1}\in\R^{d_1}$, $\bY_2,\bY_{\bP2}\in\R^{d_2}$, $d_1+d_2=d$, and with corresponding partitions of the parameters $\bxi$, $\bOmega$, $\bDelta$, and $\bDelta_\bP$. 
 Then:
\begin{itemize}
    \item[(a)] the SUN distributions of full-rank affine transformations of $\bY$ and $\bY_\bP$ coincide if and only if $\bP\in {\cal P}(m)$;
    \item[(b)] the SUN marginal distributions of $\bY_i$ and $\bY_{\bP i}$ coincide ($i=1,2$) if and only if $\bP\in {\cal P}(m)$;
    \item[(c)] the SUN conditional distributions of $(\bY_2|\bY_1=\by_1)$ and $(\bY_{\bP2}|\bY_{\bP1}=\by_1)$ coincide ($i=1,2$) if and only if $\bP\in {\cal P}(m)$.
\end{itemize}
\label{cor1}
\end{corollary}
\begin{proof} (a) If $\ba\in\R^p$ and $\bA$ is a full-rank $d\times p$ matrix, then according to \citet[p.\,199]{AC14}, $\ba+\bA^\top \bY\sim {\cal SUN}_{p,m}(\ba+\bA^\top\bxi,\bA^\top\bOmega\bA,\bDelta_\bA,\btau,\Bar{\bGamma})$ with $\bDelta_\bA=\{(\bA^\top\bOmega\bA)\odot\bI_p\}^{-1/2}\bA^\top\bomega\bDelta$; here $\odot$ denotes the element-wise product of matrices. Thus, for $\bY_\bP$,   $\bDelta_{\bA,\bP}=\{(\bA^\top\bOmega\bA)\odot\bI_p\}^{-1/2}\bA^\top\bomega\bDelta_\bP=\{(\bA^\top\bOmega\bA)\odot\bI_p\}^{-1/2}\bA^\top\bomega\bDelta\bP^\top=\bDelta_\bA\bP^\top$ and with Proposition \ref{prop1} the result is proved.\newline
(b) From \citet[eq.\,(7.6)]{AC14}, the marginal distributions of $\bY_i$ and $\bY_{\bP i}$ are  ${\cal SUN}_{d_i,m}(\bxi_i,\bOmega_{ii},\bDelta_i,\btau,\Bar{\bGamma})$  and ${\cal SUN}_{d_i,m}(\bxi_i,\bOmega_{ii},\bDelta_{\bP i},\btau_{\bP},\Bar{\bGamma}_{\bP})$, respectively ($i=1,2$). Since $\bDelta_\bP = \bDelta\bP^\top$, it follows that $\bDelta_{\bP i} = \bDelta_i\bP^\top$, and with Proposition \ref{prop1} the result is proved.\newline
(c) From \citet[eq.\,(7.7)]{AC14}, the conditional distribution of $(\bY_2|\bY_1=\by_1)$ is ${\cal SUN}_{d_2,m}(\bxi_{2\cdot 1},\bOmega_{22\cdot 1},\bDelta_{2\cdot 1},\btau_{2\cdot 1},{\bGamma}_{2\cdot 1})$ with $\bxi_{2\cdot 1}=\bxi_2+\bOmega_{21}\bOmega_{11}^{-1}(\by_1-\bxi_1)$, $\bOmega_{22\cdot 1}=\bOmega_{22}-\bOmega_{21}\bOmega_{11}^{-1}\bOmega_{12}$, $\bDelta_{2\cdot 1}=\bDelta_2-\bar\bOmega_{21}\bar\bOmega_{11}^{-1}\bDelta_1$, $\btau_{2\cdot 1}=\btau+\bDelta_1^\top\bar\bOmega_{11}^{-1}\bomega_1^{-1}(\by_1-\bxi_1)$, and ${\bGamma}_{2\cdot 1}=\bar\bGamma-\bDelta_1^\top\bar\bOmega_{11}^{-1}\bDelta_1$. Similarly, the conditional distribution of $(\bY_{\bP2}|\bY_{\bP1}=\by_1)$ is ${\cal SUN}_{d_2,m}(\bxi_{2\cdot 1},\bOmega_{22\cdot 1},\bDelta_{\bP 2\cdot 1},\btau_{\bP 2\cdot 1},{\bGamma}_{\bP 2\cdot 1})$ with:
\begin{align*}
     \bDelta_{\bP 2\cdot 1}& = \bDelta_{\bP 2}-\bar\bOmega_{21}\bar\bOmega_{11}^{-1}\bDelta_{\bP 1}=\bDelta_2\bP^\top-\bar\bOmega_{21}\bar\bOmega_{11}^{-1}\bDelta_1\bP^\top=\bDelta_{2\cdot 1}\bP^\top,\\
     \btau_{\bP 2\cdot 1} & =\btau_\bP+\bDelta_{\bP 1}^\top\bar\bOmega_{11}^{-1}\bomega_1^{-1}(\by_1-\bxi_1)=\bP\btau+\bP\bDelta_1^\top\bar\bOmega_{11}^{-1}\bomega_1^{-1}(\by_1-\bxi_1)=\bP \btau_{2\cdot 1},\\
     {\bGamma}_{\bP 2\cdot 1}& =\bar\bGamma_\bP-\bDelta_{\bP 1}^\top\bar\bOmega_{11}^{-1}\bDelta_{\bP 1}=\bP\bar\bGamma\bP^\top-\bP\bDelta_1^\top\bar\bOmega_{11}^{-1}\bDelta_1\bP^\top=\bP{\bGamma}_{2\cdot 1}\bP^\top .
\end{align*}
Note that ${\bGamma}_{2\cdot 1}$ and ${\bGamma}_{\bP 2\cdot 1}$ are no longer correlation matrices, but can be renormalized to be, which will not change the conditional density; see the corrigendum of \cite{AA06}. 
Thus, with Proposition \ref{prop1} the result is proved.
\end{proof}

The aforementioned over-parameterization problem for the SUN class is not confined there, but it also holds for the CSN and the other formulations recalled in Section~\ref{intro-basics}. In essence, the reason is that, also for these other formulations,  the components of the latent variable to which the selection mechanism is applied can be permuted in the same way as for the SUN. This fact goes on the top of the already-known over-parameterization due to the possibility of arbitrarily scaling, as recalled at the beginning of Section~\ref{intro-sun}.


\section{Further Considerations and Developments}

\subsection{Constraints on parameters to enforce identifiability}
\label{constraints}

In order to enforce identifiability of the SUN distribution, one has to fix a specific order of the $m$ latent variables, i.e., a permutation matrix $\bP\in{\cal P}(m)$.

A first approach is to order the components $\tau_i$, $i=1,\ldots,m$ of $\btau$ to follow a strictly  increasing or decreasing order. This restricts the possible values of $\btau$ and the case $\btau=\0$ is not covered, but it covers the case $\bar{\bGamma}=\bI_m$.

A second approach is to arrange $\bU_0$ by a particular order of the eigenvalues of $\bar{\bGamma}$, for example, from largest to smallest, and require that all eigenvalues are distinct. The advantage of this option it that it does not restrict the possible values of $\btau$ and the case $\btau=\0$ is covered. This also means that the case $\bar{\bGamma}=\bI_m$ cannot be covered.

From a computational point of view the first approach seems easier to implement. The second approach would need to use the spectral decomposition of $\bar{\bGamma}$.


\subsection{Sub-models of the SUN that are identifiable}

Here are some sub-models of the multivariate SUN distribution that are identifiable:
\begin{itemize}
    \item[1)] Set $\btau = \tau \boldsymbol{1}_m$ $(\tau\in \R)$, $\bDelta=\bdelta \boldsymbol{1}_m^\top$ $(\bdelta\in\R^d)$ and
$\Bar{\bGamma} = (1-\rho)\bI_m + \rho \boldsymbol{1}_m \boldsymbol{1}_m^\top$, with $\rho\in(-1,1)$ (equicorrelation). Then $\bP \btau = \bP \tau \boldsymbol{1}_m =\btau$, $\bP \bDelta^\top=\bP \boldsymbol{1}_m\bdelta^\top =\bDelta^\top$ and $\bP \bar \bGamma \bP^\top = \bar \bGamma$, since $\bP \boldsymbol{1}_m = \boldsymbol{1}_m$ and $\bP\bP^\top = \bI_m$.
Therefore the ${\cal SUN}_{d,m}(\bxi,\bOmega,\bdelta \boldsymbol{1}_m^\top,\tau \boldsymbol{1}_m,(1-\rho)\bI_m + \rho \boldsymbol{1}_m \boldsymbol{1}_m^\top)$ model is identifiable. 

This scheme can be easily extended to a formulation comprising several blocks, where each block has a parameter structure of the type indicated in the previous paragraph.


 \item[2)] Set $\btau = \alpha \boldsymbol{1}_m + \beta \bj_m$ $(\alpha,\beta \in \R, \beta \neq 0)$ where $\bj_m = (0,1,\ldots,m-1)^\top$ (or $\bj_m = (1,\ldots,m)^\top$).
In this case, $\bP\btau= \alpha \boldsymbol{1}_m + \beta \bj_m'$ with $\bj_m'=\bP\bj_m \ne \bj_m$. This means that $\btau_p \ne \btau$
for all permutations matrices $\bP \ne \bI_m$.
Note that if $\beta >0$ ($\beta <0$), this orders the elements of $\btau$ from smallest (largest) to largest (smallest), 
with the advantage that it reduces the dimensionality of $\btau$ to just $2$. Therefore the ${\cal SUN}_{d,m}(\bxi,\bOmega,\bDelta,\alpha \boldsymbol{1}_m + \beta \bj_m,\bar \bGamma)$ model is identifiable.
\item[3)] An identifiable sub-model was studied in \cite{zareifard2013non} when applying the SUN distribution to spatial statistics as ${\cal SUN}_{n,n}(\0,\omega^2 \bar {\bf C}_{\btheta},\omega \delta (1 + \delta^2)^{-1/2}\bar {\bf C}_{\btheta},\0,\bar {\bf C}_{\btheta})$, where $ \bar {\bf C}_{\btheta}$ denotes a correlation matrix of the spatial field based on a correlation function parameterized by $\btheta$. Here $\omega^2$ represents the variance and jointly controls the skewness with $\delta\in\R$. In our case, we can simply replace $\bar {\bf C}_{\btheta}$ by $\bar \bOmega$, a correlation matrix, for non-spatial modeling and employ non-zero means for the variables. The subsequent identifiable sub-model is ${\cal SUN}_{n,n}(\bxi,\omega^2 \bar\bOmega, \omega \delta (1 + \delta^2)^{-1/2}\bar\bOmega,\0,\bar\bOmega)$.   
 \item[4)] Another identifiable sub-model can be derived from the flexible subclass of the closed-skew normal (FS-CSN) distribution proposed in \cite{marquez2022flexible} for modeling spatial data{; see also \cite{2016Zareifard} for a similar model in the context of graphical modelling.} The FS-CSN subclass takes the assumption that $m=d=n$ in (\ref{csn-pdf}) and reparameterizes the CSN as $\bZ \sim {\cal CSN}_{n,n}(\bmu-b \delta \sigma \tau \bar\bSigma^{1/2} \mathbf{1}_n, \sigma^2 \tau^2 \bar\bSigma, \lambda(\sigma \tau)^{-1}\bar\bSigma^{-1/2},\0,\bI_n)$ where $\bar \bSigma^{1/2}$ is the lower triangular matrix from the Cholesky factorization of the correlation matrix $\bar\bSigma$. The reason for factoring out the scalar parameters $b$, $\delta$, $\sigma$, $\tau$, and $\lambda$ is to ensure that $\text{E}(\bZ) = \bmu$ and $\text{Var}(\bZ) = \sigma^2 \bar \bSigma$. From the formulation of the FS-CSN family follows the identifiable sub-model ${\cal SUN}_{n,n}(\bxi,\bOmega,\lambda\bOmega^{-1/2},\0,\bI_n)$. 
\item[5)] Another special case of the SUN distribution that becomes identifiable is obtained following some ideas in \cite{arellano2005fundamental}. First, consider the reparameterization
$\bOmega = \bPsi^{1/2}(\bI_d +\bLambda\bLambda\T)\bPsi^{1/2}$ and $\bomega\bDelta = \bPsi^{1/2}\bLambda$,
where $\bPsi\in\R^{d\times d}$ is positive definite  and $\bLambda\in\R^{d\times m}$. In this case, $\bY = \bxi + \bPsi^{1/2}(\bLambda\bZ +\bV)$ where $\bZ= (\bU_0|\bU_0 + \btau > \0)$, with 
$\bU_0\sim \mathcal{N}_m(\0,\bar\bGamma)$  and $\bV\sim\mathcal{N}_d(\0,\bI_d)$ being independent. Next, consider the conditions that $\bar\bGamma=\bI_m$ and $\bLambda\T\bLambda$ is diagonal, with its diagonal elements ordered, e.g., from smallest to largest.
Thus, in the numerator of  (\ref{sun-pdf}), we have easily that the covariance matrix 
$$
\bar\bGamma - \bDelta\T \bar\bOmega^{-1}\bDelta = \bI_m -\bLambda\T (\bI_d + \bLambda\bLambda)^{-1}\bLambda = (\bI_m +  \bLambda\T\bLambda)^{-1}
$$ 
is also diagonal, so that the SUN density (\ref{sun-pdf}) becomes
$$
  f(\by) = \phi_d(\by - \bxi; \bPsi + \bPsi^{1/2}\bLambda\bLambda\T\bPsi^{1/2})\:\frac{\Phi_m[(\bI_m + \bLambda\T\bLambda)^{1/2}\{\btau + \bmu(\by)\}]}{\Phi_m(\btau)},\quad \by\in\R^d,
  $$
where $\bmu(\by) = \bLambda\T(\bI_d +\bLambda\bLambda\T)^{-1}\bPsi^{-1/2}(\by - \bxi)$. Note that the order condition imposed on the diagonal elements of the diagonal matrix
$\bLambda\T\bLambda$ implies that $\bLambda_\bP\T\bLambda_\bP\ne \bLambda\T\bLambda$, for all $\bP\in\mathcal{P}(m)$, where $\bLambda_\bP=\bLambda\bP\T$.
 This condition can be omitted by assuming that the elements of $\btau$ are ordered in increasing or decreasing order. 
\item[6)] Within the Bayesian framework, \citet[Theorem 1]{durante19} has shown that the posterior distribution for the coefficients of a probit model is SUN under a Gaussian prior. Because the order of these latent variables is fixed by the prior's mean vector and covariance matrix, this posterior SUN distribution is identifiable. 
\end{itemize}
Notice that 4) can be combined with 2) to form a more general sub-model. The parameter $\btau$, set to $\0$ in 3), can also be specified as in 2) to formulate a more general model. Therefore, it appears that many more identifiable sub-models can be constructed.


\subsection{Extension to the unified skew-elliptical class}

If the normal distribution  ${\cal N}_{m+d}(\0,\bar\bOmega^*)$ in (\ref{def}) is replaced by an elliptically contoured distribution ${\cal EC}_{m+d}(\0,\bar\bOmega^*,g^{(m+d)})$, where $g^{(m+d)}$ is a generator of spherical $(m+d)$-dimensional densities, the corresponding selection mechanism leads to the unified skew-elliptical (SUE) class of distributions introduced by \cite{arellano-genton2010chjs}, which extends the SUN class. A particularly relevant member of the SUE class for applications is the unified skew-$t$ (SUT) distribution.
The SUN non-identifiability issue discussed in Section~\ref{SUN_nonidentifiability} carries on for the SUE class, since the parameters of a SUE distribution are computed using expressions similar to those of the SUN; see Proposition~2.1 and the ensuing text of \cite{arellano-genton2010chjs}. 
For the same reason, the possible remedies indicated in Subsection~\ref{constraints} for the SUN apply to the SUE class as well.


\subsection{Non-identifiability for selection distributions}

The SUN and its extension to the SUE can be obtained as special cases from the family of selection distributions discussed in \cite{ABG2006}. Specifically, they can be defined by the location-scale random vector $\bY= \bxi + \bomega \bZ$, where $\bZ = (\bU_1 | \bU_0 \in B)$, with $\bU_0$ and $\bU_1$ having an elliptically contoured joint distribution ${ \cal EC}_{m+d}(\0,\bar\bOmega^*,g^{(m+d)})$ and $B=(-\btau,{\boldsymbol\infty}) = (-\tau_1,\infty)\times\cdots\times(-\tau_m,\infty)$. For an arbitrary (absolutely) continuous joint distribution of $\bU_0$ and $\bU_1$, and a Borel set $B$ of $\R^m$, the density function of $\bY$ is $|\bomega|^{-1} f_{\bZ}\{\bomega^{-1}(\by-\bxi)\}$, $\by\in\R^d$, where $f_{\bZ}(\bz) = f_{\bU_1|\bU_0\in B}(\bz)$ is a selection density:
\begin{equation*}
    f_{\bZ}(\bz)=f_{\bU_1}(\bz)\:\frac{\Prob(\bU_0 \in B | \bU_1=\bz)}{\Prob(\bU_0 \in B)},\quad\bz\in\R^d.
\end{equation*}
Letting $B_\bA = \{\bx_\bA = \bA \bx \,|\, \bx \in B\}$, it is clear that for each matrix $\bA$ of the form $\bA = \bD \bP$, with $\bD$ an $m\times m$ positive diagonal matrix and $\bP \in \mathcal{P}(m)$, we have $\Prob( \bA\bU_0 \in B_\bA) = \Prob(\bU_0 \in B)$ and $\Prob(\bA\bU_0 \in B_\bA  | \bU_1= \bz) = \Prob(\bU_0 \in B | \bU_1= \bz)$ for every $\bz\in \R^d$. Consequently, for every matrix $\bA = \bD \bP$, we have that $f_{\bZ}(\bz) = f_{\bU_1|\bU_0\in B}(\bz) = f_{\bU_1|\bA\bU_0\in B_A}(\bz)$ for all $\bz\in\R^d$, in an obvious notation.


In this sense, the result of Proposition \ref{prop1} for the SUN class and its generalization for the SUE family are direct consequences of the symmetry of the joint distribution of $\bU_0$ and $\bU_1$, and of  noting that $B_\bA = (-\bA\btau,{\boldsymbol \infty})$, $\btau\in\R^m$, so $\Prob(\bA\bU_0 \in B_\bA) = \Prob(\bA\bU_0 > -\bA\btau)=\Prob(-\bA\bU_0 < \bA\btau) = F_{\bA\bU_0}(\bA\btau)$, since $-\bA\bU_0$ has the same distribution of $\bA\bU_0$. Similarly, $\Prob(\bA\bU_0 \in B_\bA \mid \bU_1 = \bz) =F_{\bA\bU_0 |\bU_1=\bz}(\bA\btau + \bA \bmu_{0|1}(\bz))$, where $\bmu_{0|1}(\bz)$ is the conditional location of $\bU_0 | \bU_1=\bz$.



\section{Discussion}

The identifiability condition of a statistical model is crucial when one comes to the inferential stage. Although we have described various sub-models which are identifiable, in general identifiability of the SUN family does not hold.  In the cases where this problem exists but it is not handled adequately, the inferential stage runs into major difficulties, no matter what estimation criterion or numerical optimization method is adopted. This fact is likely to explain computational difficulties or odd behaviours faced by some authors.
Potentially, an instance of this type is the non-identifiable SUN model with $\btau=\0$ and $\bar\bGamma=\bI_m$, used by \cite{gupta2012estimation} for parameter estimation using the method of weighted moments. This may explain some odd behaviour of the empirical distribution of parameter estimates in their simulation study  with $d=m=2$, where some components exhibit bi- and even tri-modality.

{We have provided some suggestions in Subsection~\ref{constraints} to make the unified skew-normal model identifiable. These `remedies' will require an extensive numerical exploration. Such a study is not suitable here for space reasons and because the present note aims at clarifying a preliminary conceptual issue, before the actual inferential problem can be tackled in a subsequent project of more operational nature.}


\section*{Acknowledgments}
The research of Kesen Wang and Marc G. Genton was supported by the
King Abdullah University of Science and Technology (KAUST).


\end{document}